\pdfoutput=1
\documentclass[a4paper,11pt,twoside]{report}

\author{Sir Michael Atiyah}
\title{Edinburgh Lectures on \\[0.5ex] Geometry, Analysis and Physics}
\date{Notes by Thomas K\"{o}ppe}

\usepackage{cmap,microtype,setspace,fancyhdr}
\usepackage{amsmath,amsthm,amscd,amssymb,bm,cancel,booktabs,ccaption}
\usepackage[vmargin={3cm,4cm},hmargin={4cm,4cm},marginparwidth=3.5cm,marginparsep=0.25cm,twoside]{geometry}

\DeclareMathOperator{\Lie}{Lie}
\DeclareMathOperator{\ev}{ev}
\DeclareMathOperator{\Map}{Map}
\DeclareMathOperator{\Cl}{Cl}


\makeatletter
\renewcommand{\section}{\@startsection{section}{1}{\z@}{-30pt}{15pt}{\bfseries}}
\renewcommand{\subsection}{\@startsection{subsection}{2}{\z@}{-15pt}{15pt}{\itshape}}
\renewcommand{\paragraph}{\@startsection{paragraph}{4}{\z@}{15pt}{-15pt}{\bfseries}}
\renewcommand{\@chapapp}{lecture series}
\renewcommand{\@makechapterhead}[1]{%
  \vspace*{50\p@}%
  {\parindent \z@ \raggedright \normalfont
    \vspace{5pt}
    \ifnum \c@secnumdepth >\m@ne
        \large\scshape \@chapapp\space \thechapter 
        \par\nobreak
        \vskip 10\p@
    \fi
    \interlinepenalty\@M
    \LARGE \scshape #1\par                         
    \vspace{5pt}
    \hrule                                        
    \nobreak
    \vskip 40\p@
  }}%
\renewcommand{\@makeschapterhead}[1]{%
  \vspace*{25\p@}%
  {\parindent \z@ \raggedright
    \normalfont
    \vspace{5pt}
    \interlinepenalty\@M
    \LARGE \scshape #1\par                         
    \vspace{5pt}
    \hrule                                        
    \nobreak
    \vskip 40\p@
  }}%
\renewenvironment{thebibliography}[1]
     {\section*{\bibname}%
      \@mkboth{\bibname}{\bibname}%
      \list{\@biblabel{\@arabic\c@enumiv}}%
           {\settowidth\labelwidth{\@biblabel{#1}}%
            \leftmargin\labelwidth
            \advance\leftmargin\labelsep
            \@openbib@code
            \usecounter{enumiv}%
            \let\p@enumiv\@empty
            \renewcommand\theenumiv{\@arabic\c@enumiv}}%
      \sloppy
      \clubpenalty4000
      \@clubpenalty \clubpenalty
      \widowpenalty4000%
      \sfcode`\.\@m}
     {\def\@noitemerr
       {\@latex@warning{Empty `thebibliography' environment}}%
      \endlist}%
\makeatother
	
\newtheoremstyle{tkthm}
  {15pt}
  {15pt}
  {\itshape}
  {}
  {\bfseries}
  {.}
  {.5em}
  {}
\newtheoremstyle{tkdef}
  {15pt}
  {15pt}
  {}
  {}
  {\bfseries}
  {.}
  {.5em}
  {}

\theoremstyle{tkthm}
\newtheorem{theorem}{Theorem}[section]
\newtheorem*{theorem*}{Theorem}

\newtheorem*{lemma*}{Lemma}
\newtheorem{proposition}[theorem]{Proposition}
\newtheorem*{proposition*}{Proposition}

\newtheorem*{corollary*}{Corollary}
\newtheorem{conjecture}[theorem]{Conjecture}
\newtheorem*{conjecture*}{Conjecture}
\theoremstyle{tkdef}
\newtheorem{definition}[theorem]{Definition}
\newtheorem*{definition*}{Definition}
\newtheorem{remark}[theorem]{Remark}
\newtheorem*{remark*}{Remark}
\newtheorem*{example}{Example}
\newtheorem*{examples}{Examples}
\newtheorem{exercise}[theorem]{Exercise}
\newtheorem{remarks}[theorem]{Remarks}

\newcommand{\ce}{\mathrel{\mathop:}=}                
\newcommand{\ec}{=\mathrel{\mathop:}}                
\newcommand{\abs}[1]{\left\lvert#1\right\rvert}      
\newcommand{\norm}[1]{{\left\lVert #1 \right\rVert}} 
\newcommand{\oton}[2][n]{#2_1, \dotsc, #2_{#1}}      
\newcommand{\eqand}{\qquad\text{and}\qquad}          

\DeclareMathOperator{\rk}{rk}

\DeclareMathOperator{\Sym}{Sym}

\DeclareMathOperator{\spec}{Spec}

\DeclareMathOperator{\id}{id}
\DeclareMathOperator{\tr}{tr}

\DeclareMathOperator{\Div}{Div}
\DeclareMathOperator{\Aut}{Aut}

\DeclareMathOperator{\Hom}{Hom}

\DeclareMathOperator{\Gr}{Gr}

\newcommand{\lsem}{[\![}
\newcommand{\rsem}{]\!]}

    \RequirePackage[dvipsnames]{xcolor}
    \RequirePackage[pdftitle={Edinburgh lectures on Geometry, Analysis and Physics}, pdfauthor={Michael Atiyah, typed by Thomas Koeppe}, pdfcreator={}, pdfproducer={}, pdfpagelayout=OneColumn, colorlinks=true, linkcolor=black, urlcolor=black, citecolor=black]{hyperref}

\begin{document}

\maketitle

\setcounter{tocdepth}{1}
\tableofcontents

\chapter*{Preface}
\addcontentsline{toc}{chapter}{Preface}

\indent 
These lecture notes are based on a set of six lectures that I gave in Edinburgh in 2008/2009 and they cover some topics in the interface between Geometry and Physics.  They involve some unsolved problems and conjectures and I hope they may stimulate readers to investigate them.
\medskip

I am very grateful to Thomas K\"oppe for writing up and polishing the lectures, turning them into intelligible text, while keeping their informal nature. This involved a substantial effort at times in competition with the demands of a Ph.D thesis.  Unusually for such lecture notes I found little to alter in them.
\medskip

\hfill Michael Atiyah

\smallskip
\hfill Edinburgh, September 2010

\newpage
\pagenumbering{arabic}

\chapter{From Euclidean 3-space to complex matrices}

\hfill{December 8 and 15, 2008}

\vfill

\section{Introduction}

We will formulate an elementary conjecture for $n$ distinct points in
$\mathbb{R}^3$, which is unsolved for $n\geq5$, and for which we have
computer evidence for $n \leq 30$. The conjecture would have been
understood 200 years ago (by Gauss). What is the future for this
conjecture?
\begin{itemize}
\item A counter-example may be found for large $n$.
\item Someone (perhaps from the audience?) gives a proof.
\item It remains a conjecture for 300 years (like Fermat).
\end{itemize}
To formulate the conjecture, we recall some basic concepts from
Euclidean and hyperbolic geometry and from Special Relativity.

\section{Euclidean geometry and projective space}

The two-dimensional sphere $S^2 = \bigl\{ (x,y,z) \in \mathbb{R}^2 :
x^2+y^2+z^2=1\bigr\}$ is ``the same as'' the complex projective line
$\mathbb{CP}^1 = \mathbb{C} \sqcup \{\infty\}$, on which we have
homogeneous coordinates $[u_1:u_2]$. Stereographic projection through
a ``north pole'' $N \in S^2$ identifies $S^2 \setminus \{N\}$ with
$\mathbb{C}$, and it extends to an identification of $S^2$ with
$\mathbb{CP}^1$ by sending $N$ to $\infty$.

\begin{exercise}\label{ex.stereo}
Suppose we have two stereographic projections from two ``north poles''
$N$ and $N'$. Show that these give a map $\mathbb{CP}^1 \to
\mathbb{CP}^1$ which is a \emph{complex linear transformation}
\[ u' = \frac{au+b}{cu+d} \text{ , where $a,b,c,d\in\mathbb{C}$ and $ad-bc\neq0$.} \]
\textit{Hint:} Start by considering stereographic projection from
$S^1$ to $\mathbb{R}$ first.
\end{exercise}

\section{From points to polynomials}

We will now associate to each set of $n$ distinct points in
$\mathbb{R}^3$ a set of $n$ complex polynomials (defined up to
scaling).

\paragraph{The case $\bm{n=2}$.} Given two points $x_1, x_2 \in
\mathbb{R}^3$ with $x_1 \neq x_2$, define
\[ f(x_1, x_2) \ce \frac{x_2-x_1}{\norm{x_2-x_1}} \in S^2 \text{ ,} \]
which gives a unit vector in the direction from $x_1$ to $x_2$. Under the
identification $S^2 \cong \mathbb{CP}^1$, $f$ associates to each pair
$(x_1,x_2)$ a point in $\mathbb{CP}^1$. Exchanging $x_1$ and $x_2$ is
just the antipodal map $x \mapsto -x$ on $S^2$.

\paragraph{The general case.} Given $n$ (ordered) points $\oton{x} \in
\mathbb{R}^3$, we obtain $n(n-1)$ points in $\mathbb{CP}^1$ by defining
\begin{equation}\label{eq.defu}
   u_{ij} \ce \frac{x_j-x_i}{\norm{x_j-x_i}} \in S^2 \cong \mathbb{CP}^1 \text{ \ for all $i\neq j$.}
\end{equation}
For each $i=1,\dotsc,n$ we define a polynomial $\beta_i \in \mathbb{C}[z]$ with roots $u_{ij}$ ($j\neq i$):
\begin{equation}\label{eq.defb}
  \beta_i(z) = \prod_{j\neq i}(z-u_{ij})
\end{equation}
The polynomials $\beta_i$ are determined by their roots up to
scaling. We make the convention that if for some $j$ we have
$u_{ij}=\infty$, then we omit the $j^\text{th}$ factor, so that
$\beta_i$ drops one degree. In fact, a more invariant picture arises
if instead we consider the associated \emph{homogeneous} polynomials
$B_i \in \mathbb{C}[Z_0,Z_1]$ given by $B_i(Z_0,Z_1) = \prod_j
\bigl(V_{ij} Z_0 - U_{ij} Z_1\bigr)$, where $[U_{ij} : V_{ij}] =
[u_{ij} : 1]$, so $\beta_i(z) = B_i(z,1)$.

We are now ready to state the simplest version of the conjecture:

\begin{conjecture}[Euclidean conjecture]\label{conj.euc}
For all sets $(\oton{x}) \subset \mathbb{R}^3$ of $n$ distinct points,
the $n$ polynomials $\beta_1(z), \dotsc, \beta_n(z)$ are linearly
independent over $\mathbb{C}$.
\end{conjecture}

\begin{remark}
The condition of linear independence of the polynomials $\beta_i$ is
independent of the choice of stereographic projection in Equation
\ref{eq.defu} by Exercise \ref{ex.stereo}.
\end{remark}

\begin{example}[$n=3$]
Suppose $x_1, x_2, x_3$ are distinct points in $\mathbb{R}^3$. They
are automatically co-planar, so that $x_1, x_2, x_3 \in \mathbb{R}^2
\subset \mathbb{R}^3$. So the points $u_{ij}$ lie in some great circle
$S^1 \subset S^2 \cong \mathbb{CP}^1$.

We can choose the north pole $N$ for the stereographic projection in
Equation \ref{eq.defu} either such that all $u_{ij}$ lie in the
equator, in which case $\abs{u_{ij}} = 1$ and $u_{ji}=-u_{ij}$, or
such that all $u_{ij}$ lie on a meridian, in which case $u_{ij}\in
\mathbb{RP}^1$ and $u_{ji} = -1\bigl/u_{ij}$.

Let us stick with the first convention, so that all $u_{ij}$ lie on
the equator and we have $\abs{u_{ij}} = 1$ and $u_{ji} =
-u_{ij}$. This defines three quadratics
\begin{eqnarray*}
  \beta_1(z) & = & (z-u_{12})(z-u_{13}) = (z-u_{12})(z-u_{13}) \\
  \beta_2(z) & = & (z-u_{21})(z-u_{23}) = (z+u_{12})(z-u_{23}) \\
  \beta_3(z) & = & (z-u_{31})(z-u_{32}) = (z+u_{13})(z+u_{23})
\end{eqnarray*}
In this case we can prove Conjecture \ref{conj.euc} in two ways:
\begin{itemize}
\item By geometric methods: Represent quadratics by lines in a plane,
      then linear dependence of the $\beta_i$ is the same as
      concurrence.
\item By algebraic methods: Compute the determinant of the
      $(3\times3)$-matrix of coefficients of the $\beta_i$ and show
      that it has non-vanishing determinant.
\end{itemize}
\end{example}

For the case $n=4$, there exists a proof using computer algebra. For
$n\geq5$, no proof is known, even for \emph{co-planar} points (i.e.\
\emph{real} polynomials). A proof will be rewarded with a bottle of
champagne or equivalent. The easiest point of departure is to consider
four points in a plane.

\section{Some physics: hyperbolic geometry}

Consider again the $2$-sphere $S^2 \subset \mathbb{R}^3$, and add a
fourth variable $t$ (for ``time''):
\begin{equation}\label{eq.mink}
  x^2 + y^2 + z^2 -R^2t^2 = 0
\end{equation}
This is the \emph{metric} of Minkowski space-time. Here $R$ is the
speed of light, and Equation \eqref{eq.mink} defines a
\emph{light cone}. Our original $2$-sphere is the base of the light
cone, the ``celestial sphere'' of an observer.

The (proper, orthochronous) Lorentz group $SO^+(3,1)$ acts on $S^2
\cong \mathbb{CP}^1$ as a group of complex projective transformations
$SL(2;\mathbb{C}) \bigl/ \pm1 = PSL(2;\mathbb{C})$.

The Euclidean version of this picture is the following: The rotation
group of $\mathbb{R}^3$, $SO(3)$, acts as $SU(2)\bigl/\pm1 \ec PSU(2)
\cong PU(2)$ on $S^2 \cong \mathbb{CP}^1$ preserving the metric given
by Equation \eqref{eq.mink} (``rigid motion''). We can also see this
as the projectivisation of the action of $SU(2)$ or $U(2)$ on
$\mathbb{C}^2$, and the projectivisation map
\[ SU(2) \twoheadrightarrow PSU(2) \cong SO(3) \]
is a \emph{double cover}. This map is the
restriction to the maximal compact subgroup of the double cover
$SL(2;\mathbb{C}) \twoheadrightarrow PSL(2;\mathbb{C}) \cong
SO^+(3,1)$.

We have two different representations of $SL(2;\mathbb{C}) \cong
\widetilde{SO^+}(3,1)$ (double cover): It acts on real $4$-dimensional
space-time $\mathbb{R}^{3,1}$ by proper, orthochronous Lorentz
transformations, and it acts on complex $2$-dimensional space
$\mathbb{C}^2$ (whose elements we call \emph{spinors}). The
fundamental link between these two representations is via
\emph{projective} spinors: A (projectivised) point in
$(\mathbb{C}^2\setminus\{0\})\bigr/\mathbb{C}^\times \cong
\mathbb{CP}^1$ corresponds to a point on the base of the light cone,
$S^2$.

Consider the hyperboloid given by $x^2+y^2+z^2-R^2t^2=-m^2$. Denote
the interior of the base of the light cone by $H_m$. The metric
induced on $H_m$ has constant negative curvature, and indeed it turns
$H_m$ into a model of hyperbolic $3$-space with curvature $-1/m^2$.

The Lorentz group $SO^+(3,1)$ acts transitively on hyperbolic
$3$-space $H^3$ by isometries, and it acts by
$SL(2;\mathbb{C})\bigl/\pm1$ on the $2$-sphere at infinity.

\section{The hyperbolic conjecture.}

Given $n$ distinct, ordered points in $H^3$, define the point $u_{ij}$
as the intersection of the oriented geodesic joining $x_i$ to $x_j$
with the $S^2$ at infinity.
We define $n$
polynomials $\oton\beta$, where $\beta_i$ has roots $u_{ij}$, as
before in Equation \eqref{eq.defb} (but note that in hyperbolic space
we no longer have a notion of ``antipodal points''). This brings us to
the second, stronger version of the conjecture:

\begin{conjecture}[Hyperbolic conjecture]\label{conj.hyp}
For all sets $(\oton{x}) \subset H^3$ of $n$ distinct points, the $n$
polynomials $\beta_1(z), \dotsc, \beta_n(z)$ are linearly independent
over $\mathbb{C}$.
\end{conjecture}

\begin{remarks}\mbox{}\\[-\baselineskip]\label{rem.hyper}
\begin{itemize}
\item There is good numerical evidence for the hyperbolic conjecture.
\item The conjecture uses only the intrinsic geometry of hyperbolic
      $3$-space, so it is invariant under the group of isometries (i.e.\
      the Lorentz group).
\item A model for $H^3$ is the open ball $B^3 \subset \mathbb{R}^3$.
      We can actually forget about the geometry of $H^3$ and just
      consider the points $\oton{x}$ to lie in $B^3 \subset
      \mathbb{R}^3$. Letting the radius of the ball $B^3$ grow (which
      is equivalent to letting the curvature of the hyperbolic space
      go to zero) exhibits the Euclidean conjecture as a limiting case
      of the hyperbolic conjecture.
\end{itemize}
\end{remarks}

\begin{remark}[The ball of radius $R$]\label{rem.radR}
As we said in Remark \ref{rem.hyper} (3), we can view the hyperbolic
conjecture as a statement about points inside the unit ball $B^3$, and
more generally inside any ball $B^3_R$ of radius $R \geq 0$ -- this
corresponds to hyperbolic space of constant curvature $-1/R^2$.
\end{remark}

We might expect that if the conjecture is false, then a
counter-example would be given by a rather special configuration of
the $n$ points $\oton{x}$. The following example treats the most
special configuration, namely the collinear one.

\begin{example}
Let $\oton{x}$ be collinear in $B^3_R$, and choose complex coordinates
on the boundary $S^2$ such that all the roots of $p_1$ are at
infinity, so that $p_1(z) = 1$. But then $p_2(z) = z$, $p_3(z) = z^2$,
\ldots, $p_n(z) = z^{n-1}$, and these are clearly linearly
independent.
\end{example}

\section{The Minkowski space conjecture}

Consider two world lines $\xi_1$, $\xi_2$ in $\mathbb{R}^{3,1}$
representing world-like motion of two ``stars''. Consider the two points
$x_1$, $x_2$ on $\xi_1$, $\xi_2$, respectively, representing events
when an ``observer'' looks up into the sky and ``sees'' the other star
on his celestial sphere, and denote the points on the respective
celestial spheres $S^2\cong\mathbb{CP}^1$ by $u_{12}$ and $u_{21}$. 
\medskip

\noindent We make this more precise and
more general:

Given $n$ moving stars (i.e.\ non-intersecting world lines) $\oton\xi$
and $n$ events $x_i \in \xi_i$, let $u_{ij}$ be the point in the
celestial sphere of $x_i$ at which the past light cone at $x_i$
intersects the world line $\xi_j$. In other words, $x_i$ ``sees''
$n-1$ other stars at points $u_{ij}$ in its own celestial
sphere. Since in (flat) Minkowski space all celestial spheres can be
identified by parallel translations, we may consider all the points
$u_{ij}$ to live in the space $\mathbb{CP}^1$.

Again we form the polynomials $\beta_i$ from the roots $u_{ij}$ as in
Equation \eqref{eq.defb} and come to the third and strongest version
of the conjecture.

\begin{conjecture}[Minkowski space conjecture]\label{conj.mink}
Let $\oton\xi \subset \mathbb{R}^{3,1}$ be $n$ non-inter\-secting world
lines in Minkowski space and $\{\oton{x}\}$ a set of $n$ (distinct)
events such that $x_i \in \xi_i$ for all $i$. Then the polynomials
$\beta_1(z), \dotsc, \beta_n(z)$ are linearly independent over
$\mathbb{C}$.
\end{conjecture}

\begin{remarks}\mbox{}\\[-\baselineskip]
\begin{enumerate}
\item Since the Lorentz group is essentially $SL(2;\mathbb{C})$, the
  Minkowski space conjecture is ``physical'', i.e.\ Lorentz-invariant.
\item If all stars emerge from a ``big bang'', i.e.\ if all world
  lines meet in a point in the past, then the Minkowski space
  conjecture reduces to the hyperbolic conjecture.
\item If stars are ``static'', the Minkowski conjecture reduces to the
  Euclidean conjecture.
\item The Minkowski conjecture is true for $n=2$ ($u_{12} \neq
  u_{21}$). There is no other evidence!
\item See \cite{At1b} and \cite{AtSut} for details.
\end{enumerate}
\end{remarks}

\paragraph{Challenge.} Prove or disprove the Minkowski space conjecture for $n=3$.

\paragraph{Remarks.}
\begin{itemize}
\item[1.] Conjecture 1.6.1 refers to world lines.  These can be interpreted as world lines of particles or ``stars" in uniform motion and this gives one version of the conjecture.  A stronger version arises if we allow all  ``physical motion" (i.e. not exceeding the velocity of light).  In \cite{At1b} I produced what purported to be an elementary counterexample for  $n= 3$.  However, on closer inspection this involves motion faster than light, so the general conjecture is still open.
\item[2.] It is even tempting to consider motion on a curved space-time background but since we now have to worry about parallel transport it is not clear how to formulate a conjecture.
\end{itemize}

\section{The normalised determinant}

We begin by recalling some basic results from linear algebra.
Consider the decomposition
\[ \mathbb{R}^2 \otimes \mathbb{R}^2 \cong \mathbb{R}^4 \cong \Sym^2
   \bigl(\mathbb{R}^2\bigr) \oplus \Lambda^2\bigl(\mathbb{R}^2\bigr)
   \cong \mathbb{R}^3 \oplus \mathbb{R}^1 \text{ .}  \]
We can view the sum on the right-hand side as the decomposition of
real $(2\times2)$-matrices into symmetric and skew-symmetric parts,
and we may think of the symmetric part $\Sym^2 \bigl(\mathbb{R}^2
\bigr)$ as a space of symmetric polynomials (of degree $2$) and of the
alternating part $\Lambda^2 \bigl( \mathbb{R}^2 \bigr)$ as the
``area'' or ``determinant''. The linear group $GL(2;\mathbb{R})$ acts
on both summands and preserves this decomposition, and it acts on the
area by multiplication by the determinant. $SL(2;\mathbb{R})$ acts
trivially on the $\mathbb{R}^1$-summand.

The complex analogue of this picture is the following: The group
$SL(2; \mathbb{C})$ acts trivially on $\Lambda^2 \bigl( \mathbb{C}^2
\bigr) \cong \mathbb{C}$ and on $\Lambda^n(\mathbb{C}^n) \cong
\mathbb{C}$. (Note: $\mathbb{C}^n \cong \Sym^{n-1} \bigl( \mathbb{C}^2
\bigr)$.) The group action preserves the standard symplectic form
on $\mathbb{C}^2$.

Now suppose we have $n$ distinct points $\oton{x}$ inside a ball of
radius $R$, and the numbers $u_{ij} \in \mathbb{CP}^1$ are defined as
in Equation \ref{eq.defu}. Lift the $u_{ij}$ to any $v_{ij} \in
\mathbb{C}^2$, i.e.\ pick a vector $v_{ij}=(z_1,z_2)$ such that
$z_1\bigl/z_2=u_{ij}$. Using the standard symplectic form, we identify
$\mathbb{C}^2$ with its dual $(\mathbb{C}^2)^\vee$, and using this identification
we consider the $v_{ij}$ as one-forms. Since $u_{ij} \neq u_{ji}$, $v_{ij} \wedge
v_{ji} \neq 0$. Now fix the constant multiplier by setting
\[ p_i = \prod_{j\neq i} v_{ij} \in \Sym^{n-1}\bigl((\mathbb{C}^2)^\vee\bigr) \cong (\mathbb{C}^n)^\vee \text{ ,} \]
and define
\begin{equation}\label{eq.defDR}
  D_R(\oton{x}) = \frac{p_1 \wedge p_2 \wedge \dotsb \wedge p_n}
  {\prod_{i<j}\bigl(v_{ij}\wedge v_{ji}\bigr)} \text{ .}
\end{equation}

\begin{remarks}
Here the numerator is an element of $\Lambda^n(\mathbb{C}^n) \cong
\mathbb{C}$, concretely given by the determinant of the $(n\times
n)$-matrix of the coefficients of the polynomials $p_i$. The
denominator is a product of elements of $\Lambda^2 \bigl( \mathbb{C}^2
\bigr) \cong \mathbb{C}$. Changing the choice of $V_{ij}$ by a factor
$\lambda_{ij}$ multiplies both numerator and denominator by the same
factor $\prod_{i\neq j} \lambda_{ij}$, so $D_R$ depends only on the
points $\oton{x}$. Permuting the points $\oton{x}$ produces the same
sign change in numerator and denominator, so $D_R$ is invariant under
permutations.
\end{remarks}

\begin{definition}[Normalised determinant]\label{def.normdet}
For $n$ distinct points $\oton{x}$ in $\mathbb{R}^3$ inside a ball of
radius $R$, we define the \emph{normalised determinant} $D_R$ to be as
in Equation \eqref{eq.defDR}. (This normalization gives $D_R=1$ for
collinear points.)
\end{definition}

\paragraph{Computation of $\bm{D_R}$.} Given $n$ distinct points
$\oton{x}$ inside a ball of radius $R$, choose for each pair $i<j$
lifts $v_{ij}$, $v_{ji}$ such that $V_{ij} \wedge V_{ji} = \omega_2$.
Write each $p_i$ in terms of the monomials $t_0^{n-1-i} t_1^i$,
where $\{t_0,t_1\}$ is a basis for $\mathbb{C}^2$ satisfying $t_0
\wedge t_1 = \omega_2$. If we denote by $P$ the $(n \times n)$-matrix
whose $(i,j)$-entry is the $j$th coefficient of $p_i$, then
$D_R(\oton{x}) = \det P$ (hence the name ``normalised determinant'').

\paragraph{Properties of the normalised determinant.}
\begin{enumerate}
\item $D_R(\oton{x})$ is invariant under the $SL(2;
      \mathbb{R})$-action (i.e.\ the isometries of $H^3_R$) on the
      points $\oton{x}$, and it is continuous in $(\oton{x})$.
\item The limit $D_\infty \ce \lim_{R\to\infty} D_R$ exists and is
      invariant under the group of Euclidean motions (translations and
      rotations of $\mathbb{R}^3$).
\item $D_R(\oton{x}) = 1$ for collinear points.
\item $D_R \to \overline{D_R}$ under reflection of $\mathbb{R}^3$ (so
      $D_R$ is real for coplanar points).
\item For $n=3$,
      \[ D_\infty = \frac12 \sum_{i=1}^3 \cos^2\bigl({\textstyle\frac{A_i}2}\bigr) \text{ ,} \]
      where $A_i$ are the angles of a triangle, varying between $1$
      for collinear and $9/8$ for equilateral configurations.
      For $n\geq4$, $D_R$ is complex-valued in general.
\item $D_\infty$ is scale-invariant: $D_\infty(\oton{\lambda x}) = D_\infty(\oton{x})$ for $\lambda>0$.
\item In the hyperbolic case, $D_R(\oton{x}) \to D_R(\oton[n-1]x)$
      as $\abs{x_n} \to R$. (This generalises to the so-called
      ``cluster decomposition'': If the points $\oton{x}$ fall into
      two ``clusters'' at great distance, then $D_R$ is approximately
      the product of the $D_R$'s of the clusters.)\label{it.cluster}
\end{enumerate}

The formalism of the normalised determinant allows us to rephrase our conjectures,
and assuming normalisation we can actually state stronger forms:
\begin{itemize}
\item The Euclidean conjecture \ref{conj.euc}. Weak form: $D_{\infty}   \neq 0$. Strong form: $\abs{D_\infty} \geq 1$ after normalisation.
\item The hyperbolic conjecture \ref{conj.hyp}. Weak form: $\abs{D_R}
  \neq 0$. Strong form $\abs{D_R} \geq 1$, after normalisation, with
  equality for collinear points.
\item We also have a new conjecture, the \emph{monotonicity
  conjecture}: $\abs{D_R}$ increases with $R$ (for fixed $\oton{x}$).
\end{itemize}

\begin{remarks}
$D_R(x) = D_{\lambda R}(\lambda x)$, so the hyperbolic conjecture is
independent of $R$. So if it is true for finite $R$, then it is true
for $R=\infty$.

The Minkowski space conjecture implies the hyperbolic conjecture:
Shrink $S^2_R$ to $S^2_{R'}$, where $R' = \abs{x_n} = \max_i
\abs{x_i}$, then apply Property \eqref{it.cluster} inductively.
\end{remarks}

The normalised determinant $D_R$ can be defined for points inside any
ellipsoid $S$, in which case we denote it by $D_S$. This is because
$S$ can be changed into a standard sphere by affine linear
transformations of $R^3$ (which preserve straight lines). We can
reduce $x^2+y^2+z^2=1$ to $x^2/a^2 + y^2/b^2 + z^2/c^2 = 1$ by choice
of $a,b,c \geq 1$.

\begin{remark}[Ellipsoid version]
The Minkowski space conjecture can be stated in terms of ellipsoids:
Suppose $S' \supseteq S$ are two ellipsoids in $\mathbb{R}^3$ containing $n$
distinct points $(\oton{x})$. Then
\[ \abs{D_{S'}(\oton{x})} \geq \abs{D_{S}(\oton{x})} \text{ .} \]
To see this, consider the situation where $S' \supseteq S$ are two
light cones. Then $\abs{D_{S'}} \geq \abs{D_S}$. A physical
interpretation is that if $S'$ is the vacuum light cone and $S$ the
light cone in a medium, then $\abs{D_\text{med}} \leq
\abs{D_\text{vac}}$.
\end{remark}

\section{Relation to analysis and physics}

\paragraph{The Dirac equation.} Let $s(x)$ be a spinor field
in $\mathbb{R}^3$. The Dirac equation in vacuum is
\[ Ds = \sum_{j=1}^3 A_j \frac{\partial s}{\partial x_j} = 0 \text{ ,} \]
where $A_j$ are $(2\times2)$-matrices, $A_j^2=-1$, $A_iA_j=-A_jA_i=A_k$
(the Pauli matrices).

\paragraph{The point monopole.} Given $(\oton{x})$, consider these as locations of $n$ Dirac monopoles
and take the Dirac equation $Ds=0$ in the background field.
We need to impose suitable singular behaviour at $\oton{x}$ and decay at infinity.

We expect an $n$-dimensional space of solutions. Examine the
asymptotic behaviour at infinity: Can we find our polynomials
$\beta_i$ in this (e.g.\ as a basis of the solutions)? Would this
imply the Euclidean conjecture?

In the hyperbolic case, the asymptotic behaviour may be exponential
decay, with polynomial angular dependence. Would this imply the
radius-$R$ conjecture for finite $R$?

\paragraph{The four-dimensional variant.} Let $M^4$ be the Hawking-Gibbons
$4$-manifold, which has an action of $U(1)$. The quotient is $M^4
\bigl/ U(1) = \mathbb{R}^3$, and the $U(1)$-action has $n$ fixed
points, which determine $n$ points $\oton{x}$ in $\mathbb{R}^3$.

Reinterpret on $M^4$: The solutions of the four-dimensional Dirac
equation on $M^4$ inherit an action of $U(1)$. The invariant solutions
on $M^4$ correspond to the singular solutions on $\mathbb{R}^3$. This
disposes of the singular behaviour at $x_i$. We still require decay at
infinity.

Next step: The Dirac equation is conformally invariant, so we can form
the conformal compactification $\overline{M}$ (which has a mild
singularity at infinity). This replaces asymptotic behaviour by local
behaviour near infinity.

In the final step, we form the twistor space of $\overline{M}$ and use
the complex methods of sheaf theory: Under the twistor transform,
solutions of the Dirac equation correspond to sheaf cohomology. In
particular, we expect a certain first cohomology to have dimension
$n$.

To relate this to polynomials and our conjectures, we must use
\emph{real} numbers and positivity. This is close to (real) algebraic
geometry.

\paragraph{Hyperbolic analogue.}
Four-manifold $N^4$ with special metric and $U(1)$-action with $n$
fixed points and $N^4 \bigl/ U(1) = B^3$, the ``inside'' of $S^2$ in
$\mathbb{R}^3$ with the hyperbolic metric. It admits a conformal
compactification $\overline{N}$, on which we have a $U(1)$-action with
$n$ fixed points and a fixed $S^2 \subset \overline{N}\setminus
N^4$. Twistor methods still apply to this case, but is it better than
the Euclidean case? This leads to the theory of LeBrun manifolds. See
Atiyah-Witten, which includes a problem about the existence of $G_2$
metrics on $7$-manifolds which are $\mathbb{R}^3$-bundles over $N^4$
and generalise the cases $n=0$ and $n=1$.

\paragraph{Lie group generalisation.}
The Euclidean conjecture implies the existence of a continuous map
\[ f_n \colon C_n\bigl(\mathbb{R}^3\bigr) \to GL(n; \mathbb{C}) \bigl/
   (\mathbb{C}^\times)^n \to U(n )\bigl/ T^n \]
compatible with the action of the symmetric group. Specifically, the
value in $GL(n; \mathbb{C})$ is the matrix of coefficients of the
polynomials $p_i$, and the quotient by $(\mathbb{C}^\times)^n$ accounts
for the freedom of scale.

The configuration space can be described as follows.
\[ C_n\bigl( \mathbb{R}^3 \bigr) = \Lie\bigl(T^n\bigr) \otimes \mathbb{R}^3 \setminus \mathcal{S} \text{ ,} \]
where the $\mathbb{R}^3$-factor contains the coordinates of the points, the factor $\Lie\bigl(T^n\bigr)$
accounts for the $n$ points, and $\mathcal{S}$ is the union of codimension-$3$ linear subspaces
$\mathcal{S}_\alpha$, where $\mathcal{S}_\alpha$ is the kernel of the linear map (the root map)
\[ \alpha \otimes \id_{\mathbb{R}^3} \colon \Lie\bigl(T^n\bigr) \otimes \mathbb{R}^3 \to \mathbb{R}^3 \]
extending the roots $\alpha$ of $U(n)$, accounting for the fact that
the $n$ points are required to be \emph{distinct}. (The roots of
$U(n)$ are formed by elements $x_i-x_j$. Note that
$\Lie\bigl(T^n\bigr)$ is the Cartan subalgebra of $\mathfrak{u}(n)$.)

This leads us to a generalisation of our conjectures. Let $G$ be a
compact Lie group (e.g.\ $SO(n;\mathbb{R})$) and $G_\mathbb{C}$ its
complexification (e.g.\ $SO(n;\mathbb{C})$). Let $T \leq G$ be a
maximal torus with complexification $T_\mathbb{C}$, and let $W \ce
N(T)\bigl/T$ be the Weyl group of $G$, which permutes the roots.

\begin{conjecture}[Lie group conjecture]\label{conj.lie}
If $G$ is a Lie group as above with rank $n$, then there exists a
continuous map
\[ f_n \colon \Lie(T) \otimes \mathbb{R}^3 = \mathcal{S} \to
   \frac{G_\mathbb{C}}{T_\mathbb{C}} \to \frac{G}{T} \]
compatible with the action of the Weyl group $W$.
\end{conjecture}

In a joint paper with Roger Bielawski (\cite{AtBi}) we used Nahm's equations
\[ \frac{dA_1}{dt} = [A_2,A_3] \text{ (and cyclic permutations),} \]
where $A_i \colon (0,\infty) \to \Lie(G)$ are functions of $t$ subject
to suitable boundary conditions $t\to0$, $t\to\infty$, to prove the
existence of a map to $G \bigl/ T$. Problems:
\begin{enumerate}
\item For $G=U(n)$, is this the same as a map given by polynomials?
\item Is there an explicit algebraic analogue for a map to
      $G_\mathbb{C} \bigl/ T_\mathbb{C}$?.
\item Is there any generalisation of the hyperbolic conjecture from
      $GL(n;\mathbb{C})$ to other Lie groups?
\end{enumerate}

\section{Mysterious links with physics}
\begin{itemize}
\item Origin in Berry-Robbins on spin statistics.
\item Link to Dirac equation?
\item Generalisation to Minkowski space.
\item Nahm's equations and gauge theory.
\item Link to Hawking-Gibbons metric?
\item Twistor interpretation?
\end{itemize}

\paragraph{Key fact of physics.} The base of the light cone is $\mathbb{CP}^1$.
It is Penrose's philosophy that this must be the origin of complex
numbers in quantum theory, and it must lie behind any unification of
General Relativity and Quantum Mechanics.

\bigskip\noindent
What is the physical meaning of our conjectures?

\section*{List of conjectures}
\addcontentsline{toc}{section}{List of conjectures}
\begin{itemize}
\item Conjecture \ref{conj.euc}: The Euclidean conjecture (weak and strong).
\item Conjecture \ref{conj.hyp}: The hyperbolic conjecture (weak and strong).
\item The monotonicity conjecture for the normalised determinant.
\item Conjecture \ref{conj.mink}: The Minkowski space conjecture.
\item Conjecture \ref{conj.lie}: The Lie group conjecture.
\end{itemize}

\chapter{Vector bundles over algebraic curves and counting rational points}

\hfill{February 9, 16, 23 and March 2, 2009}

\vfill

\section{Introduction}

There are two themes, both initiated by A.\ Weil:
\begin{enumerate}
\item Extension of classical ideas in algebraic geometry, number
      theory, physics from Abelian (scalars, $U(1)$) to non-Abelian
      (matrices, $U(n)$) settings.
\item Connection between homology and counting rational points over finite fields.
\end{enumerate}

\section{Review of classical theory}

(Abel, Jacobi, Riemann, \ldots) Consider complex projective space
\[ \mathbb{CP}^{n-1} \equiv \mathbb{P}^{n-1} \ce (\mathbb{C}^n\setminus \{0\}) \bigl/ \mathbb{C}^\times \]
with homogeneous coordinates $[z_1: \dotsc: z_n]$. Rational functions
on $\mathbb{P}^{n-1}$ are fractions $f(\oton{z})\bigl/g(\oton{z})$,
where $f$ and $g$ are homogeneous polynomials of the same degree.

Note that a meromorphic function on $\mathbb{P}^{n-1}$ is determined
up to scale by its zeros and poles (Liouville). On projective space,
global complex analysis is just algebraic geometry (Serre).

There is the standard line bundle $L$ over $\mathbb{P}^{n-1}$, i.e.\
$L \cong \mathcal{O}_{\mathbb{P}^{n-1}}(1)$. Holomorphic sections of
its $k^\text{th}$ power $L^k = L \otimes L \otimes \dotsb \otimes L$
are just homogeneous polynomials of degree $k$.

\paragraph{Algebraic curves.} Let $n=3$, so we consider the projective
plane $\mathbb{P}^2$. A curve of degree $k$ is given as the locus of
points $z$ such that $f(z_1, z_2, z_3) = 0$, where $f$ is a
homogeneous polynomial of degree $k$. Non-singular curves are just
compact Riemann surfaces, so topologically they are entirely
determined by its genus $g$. A Riemann surface $X$ of genus $g$ has
first Betti number $b_1 = \dim H_1(X;\mathbb{Q}) = 2g$.

If $X$ is a curve of degree $k$ with double points, then removing
these double points leaves a Riemann surface. We have a formula
\[ g = \frac12(k-1)(k-2) - \delta \text{ ,} \]
where $\delta$ is the number of double points. If $\delta=0$, then for
$k=1,2$ we find that $X$ is a rational curve, i.e.\ $g=0$; and for
$k=3$ we get an elliptic curve with genus $g=1$.

Another interpretation is that $g$ is the dimension of the space of
holomorphic differentials (which look locally like $\phi(z)\,d\!z$,
where $\phi$ is holomorphic). When $g=0$, the curve is the Riemann
sphere $\mathbb{CP}^1 = \mathbb{C} \cup \{\infty\}$, and the
differential $d\!z$ has a pole at infinity, so it is not
holomorphic. When $g=1$, the curve is the torus
$\mathbb{C}\bigl/\mathbb{Z}^2$, so the differential $d\!z$ on
$\mathbb{C}$ descends to a holomorphic differential on $X$.

\paragraph{Period matrices.} Let $\oton[g]\omega \in H^1(X; \mathbb{C})$
be a basis of holomorphic differentials and $\oton[2g]\alpha \in
H_1(X; \mathbb{Z})$ a basis for the $1$-cycles of a genus-$g$ curve
$X$. The $(g\times 2g)$-matrix with entries $\int_{\alpha_j}\omega_i$
is called the \emph{period matrix} of $X$.

\paragraph{Divisors.} We call a subvariety of codimension $1$ a
\emph{divisor}. Since curves are $1$-dimensional, divisors on curves
are just points. The free Abelian group of all divisors of a variety
$X$ is denoted by $\Div(X)$, and so if $X$ is a curve, elements of
$\Div(X)$ are just formal sums $D = \sum_{i=1}^{N} n_i P_i$, where
$P_i \in X$ are points. The \emph{degree} of such a divisor $D$ on a
curve is defined as $\deg D \ce \sum_{i=1}^N n_i$.

\paragraph{Jacobians.} The \emph{Jacobian} of $X$, written $J(X)$, is a
complex torus of complex dimension $g$, given as
\[ J(X) = \mathbb{C}^g\bigl/\text{lattice} = \text{hol.\ differentials}
  \bigl/ \text{differentials with integer periods} \text{ .} \] The
significance of the Jacobian lies in the following observation.  Let
$\phi$ be a rational (meromorphic) function on $X$, and define the
divisors
\begin{eqnarray*}
 D_0(\phi) &\ce& \text{set of zeros of $\phi$, with multiplicities,} \\
 D_\infty(\phi) &\ce& \text{set of poles of $\phi$, with multiplicities,} \\
 D(\phi) &\ce& D_0(\phi) - D_\infty(\phi) \text{ \ (the \emph{divisor of $\phi$}).}
\end{eqnarray*}
Then $\deg D_0(\phi) = \deg D_\infty(\phi)$. This motivates the
question for the converse: Given two divisors $D_1$ and $D_2$ of the
same degree, when does there exist a function $\phi$ on $X$ with
$D_0(\phi) = D_1$ and $D_\infty(\phi) = D_2$?

This is always true for $g=0$, but not otherwise. The ``gap'' between
divisors of degree zero and divisors of meromorphic functions is
measured precisely by the \emph{divisor class group} $\Cl(X)$. The degree-$0$
part of it is
\[ \Cl^0(X) \ce \frac{\text{divisors of degree $0$}}
   {\text{{divisors of functions}}} \text{ .} \]
(Divisors of the form $D=D(\phi)$ are also called \emph{principal
divisors}.) For $g=0$, the group $\Div^0(X)$ is trivial, but for $g=1$, the
divisor class group is precisely the Jacobian (or its dual) -- this is
the content of the Abel-Jacobi Theorem. Moreover, the group
$\Cl^0(X)$ is the group of isomorphism classes of holomorphic line
bundles of degree $0$, which are just given by elements of
\[ \Hom\bigl(\pi_1(X), U(1)\bigr) \]
(up to duality and complex structure). Differential geometry shows
that a holomorphic line bundle of degree zero (i.e.\ first Chern class
zero) has a unique flat unitary connection.

\bigskip\noindent This is the beginning of the link with
physics. Maxwell's equations deal with the curvature of a line bundle
on space-time.

\section{Analogy with number theory}

\noindent\hfil\begin{tabular}{rl}\toprule
\emph{Number theory} & \emph{Algebraic geometry} \\\midrule
Ring of integers $\mathbb{Z}$  & complex (affine) line, the ring $\mathbb{C}[z]$ \\
primes & points \\
factorisation of integers & factorisation of polynomials \\
``infinite prime'' & point at infinity in $\mathbb{P}^1$ \\
algebraic number field & algebraic curve (covering of a line) \\
lack of unique factorisation & not all divisors come from functions \\
ideal class group & divisor class group \\
Galois group & $\pi_1(X)$ \\\bottomrule
\end{tabular}\hfil

\bigskip\noindent
The classical analogy is that between the ring of integers in number
theory and the polynomial rings in geometry. A ``half-way house'' is
an \emph{algebraic curve over a finite field}. A finite field is a
field $\mathbb{F}_q$ with $q$ elements, where $q=p^n$ for some prime $p$.

We have ``function field analogues'' of geometric statements, e.g.\ a
Riemann hypothesis (which is proved for finite fields; also for
algebraic varieties of any dimension).

The key fact for algebraic geometry over $\mathbb{F}_q$ is the
existence of the \emph{Frobenius map} $x \mapsto x^q$. (Recall that in
characteristic $p$, $(x+y)^p = x^p + y^p$.) There is no such analogue
in characteristic zero (but physics suggest rescaling the
metric\footnote{A quick explanation of this remark is in order: In
  differential geometry, a differential form scales with the power of
  its degree, so rescaling picks out the degree of the form. In
  characteristic $p$, the eigenvalues of the Frobenius map pick out
  the dimension of the cohomology.}).

\section{Relation between homology and counting rational points}

\begin{definition}[Poincar\'e series]
For any topological space $X$ whose singular homology groups
$H_k\bigl(X;\mathbb{Q}\bigr)$ are finite-dimensional vector spaces, we
define the \emph{Poincar\'e series} of $X$ to be the formal power series
\[ P_X(t) \ce \sum_{k=0}^\infty \dim H_k(X;\mathbb{Q}) \; t^k \text{ .} \]
\end{definition}

\begin{proposition}
If $X$ is a manifold or homotopy-equivalent to a manifold, then $P_X$
is in fact a polynomial.
\end{proposition}

\begin{example}
Consider the space $\mathbb{P}^{n-1}$. Over $\mathbb{C}$, this has
Poincar\'e series
\begin{equation}\label{eq.poin}
  P(t) = 1 + t^2 + t^4 +\dotsb + t^{2n-2} \text{ .}
\end{equation}
Over $\mathbb{F}_q$, the number of points in $\mathbb{P}\mathbb{F}_q^{n-1}$ is
\[ \frac{q^n-1}{q-1} = 1 + q + q^2 + \dotsb + q^{n-1} \text{ .} \]
This agrees with \eqref{eq.poin} if we put $q=t^2$. Note that since we
can replace $q$ by $\tilde{q} = q^n$, $n=1,2,\ldots$, we can think of
$q$ as a variable like $t$. This extends to all algebraic varieties.
\end{example}

\begin{exercise}
Check that a similar relation between the Poincar\'e series over
$\mathbb{C}$ and the number of points over a finite field
$\mathbb{F}_q$ holds for the full flag variety $U(n) \bigl/ T^n$.
(Hint: Use successive fibrations by projective spaces.)
\end{exercise}

\paragraph{Generalisation from $U(1)$ to $U(n)$.} This
corresponds to generalising from line bundles to vector bundles. In
number theory, this corresponds to non-Abelian class field
theory. There are representations from the Galois group to $U(n)$,
Langlands programme\ldots In physics, this is related to non-Abelian
gauge theories and Yang-Mills theory.

Returning to algebraic geometry, we will focus on an algebraic curve
$X$ (either over $\mathbb{C}$ or over $\mathbb{F}_q$). The Jacobian is
replaced by a ``moduli space'' of vector bundles over $X$. There are a
few difficulties:
\begin{itemize}
\item There is no group structure (the tensor product does not
      preserve rank for ranks $>1$).
\item Bundles of rank $n$ can decompose into bundles of lower rank.
\end{itemize}

There is a moduli space $M_s(X,n,k)$ of holomorphic rank-$n$ bundles
of degree $k$ which are \emph{stable}. Here $k$ is the degree of the
determinant line bundle, which is the first Chern class; in symbols:
$\deg E \ce \deg \Lambda^n E \equiv c_1(\Lambda^n E)$. The space
$M_s(X,n,k)$ is a compact algebraic variety if $\gcd(n,k)=1$, e.g.\ if
$n=2$, $k=1$.

For $k=0$, the space $M_s(X,n,0)$ is the space of irreducible
representations $\pi_1(X) \to U(n)$. To see this, note that such a
representation is a choice of $2g$ unitary matrices $\oton[g]A,
\oton[g]B \in U(n)$ such that $\prod_{i=1}^g [A_i,B_i] = 1$, modulo
conjugation by $U(n)$.

For a general $k$, we replace this condition by $\prod_{i=1}^g [A_i, B_i] =
\zeta \id$, where $\zeta$ is a central element of $U(n)$ and $\zeta^k=1$.
(For example, for $n=2$, $k=1$, we have $\prod_i [A_i, B_i] = -\id$.)

The general problem is to study $M_s$.
\begin{enumerate}
\item What does $M_s$ look like topologically?
\item What are its Betti numbers $b_i = \dim H_i\bigl(M_s;\mathbb{Q}\bigr)$?
\end{enumerate}

For a connected, oriented manifold $X$, we have the Poincar\'e
polynomial $P_X(t) = \sum_{i=0}^{\dim X} = b_i t^i$. Note that for any
compact manifold $X$, we have $\deg P_X = \dim X$, and $P_X$ is
palindromic (by Poincar\'e duality). Furthermore, $P_{X \times Y}(t) =
P_X(t)P_Y(t)$.

\begin{example}
If $X$ is a Riemann surface of genus $g$, then $P_X(t) = 1+2\,g\,t +
t^2$, and $P_{J(X)}(t) = (1+t)^{2g}$.
\end{example}

What is $P_{M_s(X,n,k)}(t)$? What is it when $\gcd(n,k)=1$? Let us
consider the special case $n=2$, $k=1$ on a curve $X$ of genus
$g(x)=2$. Then
\[ P_{M_s(X,2,1)}(t) = \bigl(1+t^2+4t^3+t^4+t^6\bigr) \bigl(1+t\bigr)^4
   = \bigl(1+t^2+4t^3+t^4+t^6\bigr) P_{J(X)}(t) \text{ .} \]

\paragraph{Note.} For $n=2$ and $g\geq2$, we have $\dim_\mathbb{C}
M_s(X,2,k) = (3g-3)+g$, and $g = \dim J(X)$.

\bigskip\noindent Let $\det \colon M_s(X,n,0) \to J(X)$ be the
determinant map $E \mapsto \det E \equiv \Lambda^n E$, and denote by
$M^0$ the fibre of $\det$ over some point. We have a general result.

\begin{theorem}[Formula for general $g \geq 2$]\label{thm.pms}
\begin{equation}\label{eq.pms}
  P_{M^0}(t) = \frac{(1+t^3)^{2g}}{(1-t^2)(1-t^4)} - \frac{t^{2g}(1+t)^{2g}}{(1-t^2)(1-t^4)}
\end{equation}
\end{theorem}
\begin{exercise}
This should be a palindromic polynomial of degree $6g-6$, all of whose
coefficients are non-negative. Prove this.
\end{exercise}

Let us write in short $M_g(n,k)$ for $M_s(X,n,k)$, the moduli space of
stable vector bundles of rank $n$ and degree $k$ on a smooth curve $X$
of genus $g$. What can we say for $n \geq 2$? For $\gcd(n,k)=1$,
$M_g(n,k)$ is a complex manifold of dimension
$(3g-3)+g$. Topologically, $M_g(n,k)$ is given by $\oton[g]A,\oton[g]B
\in U(n)$ such that $\prod_{i=1}^g [A_i,B_i] = \sigma$, where $\sigma
= e^{2\pi i/n}$, modulo conjugation by $U(n)$.

\paragraph{Specific question.} What is the homology of $M_g(n,k)$?
What is its Poincar\'e polynomial? Recall:
\[ P_M(t) \ce \sum_{i=0}^N \dim H^i\bigl(M_g(n,k);\;\mathbb{Q}\bigr) \; t^i \]
Here $N = 8g-6$. For $n=1$, $P_{M_g(1,k)}(t) = (1+t)^{2g}$, independent of $k$.

For $n=2$, $k=1$, the moduli space decomposes as $M_g(2,1) =
M^0_g(2,1) \times J(X)$, and the Poincar\'e polynomial of $M^0_g(2,1)$
is given by Equation \eqref{eq.pms}.

\section{The approach via Morse theory}

\subsection{Basic Morse theory}

Let $Y$ be an $n$-dimensional manifold and $f \colon Y \to \mathbb{R}$
a function; the points $x \in Y$ where $df(x) = 0$ are called the
\emph{critical points of $Y$}. The Hessian, which we write
briefly as ``$d^2\!f$'', is a quadratic form, and we call $f$ a
\emph{Morse function} if $d^2\!f$ is non-degenerate at all critical
points of $f$. By the Morse Lemma, there exist near every critical
point $p$ local coordinates $\{x_i\}$ in which $f$ takes the form
\[ f(p+x) = f(p) -x_1^2 - x_2^2 -\dotsb - x_r^2 + x_{r+1}^2 + \dotsb + x_n^2 \text{ .} \]
The integer $r$ is called the \emph{Morse index} of the critical
point.  If $r=0$, $f$ has a minimum; if $r=n$, $f$ has a maximum, and
if $0<r<n$, $f$ has a saddle point.

If $f$ is a Morse function on $Y$, the \emph{Morse polynomial} is
\[ M_{Y,f}(t) = \sum_{Q} t^{\gamma(Q)} \text{ ,} \]
the sum over all non-degenerate critical points $Q$, and $\gamma(Q)$
is the Morse index of $Q$. It can be shown that
\[ M_{Y,f}(t) \geq P_Y(t) \text{ ,} \]
with equality in ``good cases''.

\begin{examples}\mbox{}\\[-\baselineskip]
\begin{itemize}
\item Let $Y = S^1$ and $f \colon Y \to \mathbb{R}$ the height
      function. Then $M_{Y,f}(t) = P_Y(t) = 1+t$; this is a ``good case''.
\item Let $Y = S^1$, but ``pinched'', and $f$ again the height
      function. Then $M_{Y,f}(t) = 2+2t$, a ``bad case''.
\item Let $Y = S^1 \times S^1$ be the torus and $f$ the height
      function. Then $M_{Y,f}(t) = 1+2t+t^2 = (1+t)^2 = P_Y(t)$,
      another ``good case''.
\item Let $Y = \mathbb{CP}^{n-1}$ and
      \[ f(z) = \frac{\sum_{i=1}^n \lambda_i\abs{z_i}^2}{\sum_{i=1}^n
         \abs{z_i}^2} \text{ \ with \ } \lambda_1 < \dotsb < \lambda_n \text{ .} \]
      Then the critical points of $f$ are $Q_j$ where $z_j=1$ and
      $z_i=0$ for $i \neq j$, with indices $\gamma(Q_j) = 2j-2$. Hence
      $M_{Y,f}(t) = 1 + t^2 + \dotsb + t^{2n-2} = P_Y(t)$, and we have
      another ``good case''.
\end{itemize}
\end{examples}

We generalise the notion of non-degeneracy to allow critical
submanifolds. $Q \subseteq Y$ is a critical submanifold if $df=0$
along $Q$ and $d^2f$ is non-degenerate in normal directions. The Morse
index of $Q$, written again as $\gamma(Q)$, is the number of linearly
independent negative normal directions. Such a function will be called
a \emph{Morse-Bott function}.

\begin{definition}
If $f \colon Y \to \mathbb{R}$ is a Morse-Bott function, the
\emph{Morse polynomial} of $f$ is
\[ M_{Y,f}(t) = \sum_Q t^{\gamma(Q)} P_Q(t) \text{ ,} \]
where the sum is taken over all non-degenerate critical submanifolds
$Q \subset Y$.
\end{definition}

Again we have the Morse inequality $M_{Y,f}(t) \geq P_Y(t)$, with
equality in good cases.

\begin{examples}\mbox{}\\[-\baselineskip]
\begin{itemize}
\item Let $Y = \mathbb{CP}^{n-1}$ and
      \[ f(z) = \frac{\sum_{i=1}^n \lambda_i\abs{z_i}^2}{\sum_{i=1}^n
         \abs{z_i}^2} \text{ \ with \ } \lambda_1 \leq \dotsb \leq
         \lambda_n \text{ , $\lambda_1\neq\lambda_n$.} \]

      If for example $\lambda_1=\lambda_2=\dotsb=\lambda_{n-1} <
      \lambda_n$, then $Q_\text{min} = \mathbb{CP}^{n-2}$ and
      $Q_\text{max} = \{\text{pt.}\} = [0:\dotsc:0:1]$, and so
      \[ M_{Y,f}(t) = P_{\mathbb{CP}^{n-2}}(t) + t^{2n-2} =  P_Y(t) \text{ .} \]

      More generally, if $\lambda_1 = \dotsb = \lambda_r <
      \lambda_{r+1} < \dotsb \lambda_n$, then
      \[ M_{Y,f}(t) = P_{\mathbb{CP}^{r-2}}(t) + t^{2r} + \dotsb +  t^{2n-2} =  P_Y(t) \text{ .} \]
\item Now take $n=\infty$ in the last example. Then
      $P_{\mathbb{CP}^\infty}(t) = 1+t^2+\dotsb =
      \frac{1}{1-t^2}$. But we still have
      \[ P_{\mathbb{CP}^\infty}(t) = P_{\mathbb{CP}^{r-1}}(t) + t^{2r} + t^{2r+2} + \dotsb \text{ ,} \]
      so we conclude that $P_{\mathbb{CP}^{r-1}}(t) = \frac{1}{1-t^2} - \sum_{k=r}^\infty t^{2r}$.
\end{itemize}
\end{examples}

The last example is the prototype of the method to compute
$P_{Q_\text{min}}(t)$ of some critical manifold $Q_\text{min}$ in
terms of the (possibly infinite-dimensional) total space and higher
critical points. We will use this method again later to compute the
Poincar\'e series of the moduli space of $U(2)$-bundles over a curve
of genus $g$.

\subsection{Equivariant cohomology, or The effect of symmetry}

Let $G$ be a compact Lie group (for instance $U(1)$ or $U(n)$) and
suppose $G$ acts on a manifold $Y$. If the action is free, then $Y/G$
is a manifold and has nice cohomology and Poincar\'e series. If the
action is not free, $Y/G$ has singularities. What to do?

\begin{definition}[Equivariant cohomology]
We define $H_G^*(Y) \ce H^*(Y_G)$ to be the \emph{$G$-equivariant
cohomology of $Y$}, where $Y_G$ is given by the \emph{Borel construction}
\[ Y_G \ce (EG \times Y) \bigl/ G \text{ ,} \]
where $EG$ is a contractible space with a free $G$-action, and the
action of $G$ on $EG\times Y$ is $g.(e,y) = (g.e, g.y)$. (In fact,
$EG$ is the total space of the classifying fibration $G
\hookrightarrow EG \twoheadrightarrow BG$.)
\end{definition}

\begin{example}
Let $G=U(1)$ and $EG = \mathbb{C}^\infty\setminus\{0\} = \varinjlim
\bigl(\mathbb{C}^N \setminus \{0\} \bigr)$. Then $BG \ce EG \bigl/ G =
\mathbb{CP}^\infty$. We compute:
\[ H_G^*(\text{pt.}) = H^*\bigl(\mathbb{CP}^\infty\bigr) \eqand
   P(t) = 1+t^2+\dotsb = \frac{1}{1-t^2} \text{ .} \]
\end{example}

Note that the projection
\[ Y_G = (EG \times Y) \bigl/ G \to EG\bigl/G \ec BG \simeq \{\text{pt.}\}_G \]
gives a homomorphism
\[ H^*_G\bigl(\{\text{pt.}\}\bigr) = H^*\bigl(BG\bigr) \longrightarrow H^*_G\bigl(Y\bigr) \text{ ,} \]
which turns $H^*_G\bigl(Y\bigr)$ into a graded module over the graded
cohomology ring $H^*_G(\text{pt.})$. We saw from the example that for
$G = U(1)$, the equivariant cohomology $H^*_G(\text{pt.}) =
H^*\bigl(\mathbb{CP}^\infty\bigr)$ is a polynomial ring in one
variable $u$ of degree $2$, and we may take $u$ to be the Chern class
of the tautological line bundle on $\mathbb{CP}^\infty$. More
generally, for $G = U(n)$ the equivariant cohomology
$H^*_G(\text{pt.}) = H^*\bigl(BU(n)\bigr)$ is a polynomial ring in $n$
variables $\oton{u}$ of degrees $2,4,\dotsc,2n$, and again the $u_i$ may
be interpreted as the Chern classes of the tautological $n$-plane
bundle over $BU(n) = \Gr_n\bigl(\mathbb{C}^\infty\bigr)$.

\begin{definition}
Let $Y$ be a manifold with an action of a compact Lie group $G$ as
above. The \emph{equivariant Poincar\'e series} of $Y$ is
\[ P_Y^G(t) = \sum_{k=0}^\infty \dim H^k_G\bigl(Y\bigr) \, t^k \text{ .} \]
\end{definition}

\begin{remark}
If the action of $G$ on $Y$ is free, then $Y_G \cong (EG \times Y)
\bigl/ G \simeq Y\bigl/G$, and so $P^G_Y(t) = P_{Y/G}(t)$ is a
polynomial. In general, however, $P^G_Y(t)$ is only a power series
which is the expansion of a rational function. If $Y$ is
contractible, then $Y_G \simeq BG$, and so $H^*_G\bigl(Y\bigr) =
H^*\bigl(BG\bigr)$.
\end{remark}

\paragraph{Equivariant Morse theory.}
Suppose $G$ acts on $Y$ and $f \colon Y \to \mathbb{R}$ is a
$G$-invariant Morse-Bott function, i.e.\ $f$ is a Morse-Bott function
and $f(g.y) = f(y)$ for all $g \in G$. If $G$ acts freely on $Y$, then
$f$ induces a function $f_G \colon Y\bigl/G \to \mathbb{R}$, and we
can apply Morse theory to $f_G$. Otherwise, consider $f$ on $Y$, but
remember the $G$-action and use $H_G$, that is, consider $f$ as a
Morse function on $Y_G$.

\begin{example}
Let $Y=S^2$ and $G=U(1)$, acting by a simple rotation with two fixed
points, and let $f$ be the height function. Then
\[ M_{Y,f}^G(t) = \underbrace{\frac{1}{1-t^2}}_\text{min.} + \underbrace{
   \frac{t^2}{1-t^2}}_\text{max} = \frac{1+t^2}{1-t^2} \text{ .} \]
This is a ``good case'', since we also have $P^G_Y(t) = (1+t^2)\bigl/(1-t^2)$.
\end{example}

\paragraph{Some criteria for a good Morse(-Bott) function.} The following
conditions allow us to conclude that a Morse polynomial (or power
series) is ``good'', i.e.\ equal to the Poincar\'e series.
\begin{itemize}
\item If all Morse indices and all Betti numbers are even. (E.g.\ for $\mathbb{CP}^{n-1}$.)
\item In the equivariant case: If each critical submanifold is
      point-wise fixed by a some $U(1) \subset G$ which has no fixed
      vectors in the negative normal bundle.
\end{itemize}
We will use these criteria in gauge-theoretical computations in the
following section.

\subsection{Application to infinite dimensions (gauge theory)}

Let $X$ be a surface of genus $g \geq 2$ and $A$ a $G$-connection for a
vector bundle of rank $n$ over $X$, where $G = U(n)$. For the trivial
bundle $X \times \mathbb{C}^n$,
\[A = \sum_{i=1}^2 \, A_i(x) \, dx_i \text{ ,} \]
where $(x_1, x_2)$ are local coordinates on $X$ and $A_i \in
\mathfrak{u}(n)$, the Lie algebra of skew-Hermitian $(n\times
n)$-matrices. The \emph{curvature} of the connection is (locally, or
globally in the case of the trivial bundle)
\[ F_A = dA + A \wedge A \in \Omega^2\bigl(X; \mathfrak{u}(n)\bigr) \text{ .} \]
The Lie algebra $\mathfrak{u}(n)$ admits an invariant inner product,
so we can define a norm $\norm{-}$ on it. The \emph{Yang-Mills
functional} of the connection $A$ is
\[ \phi(A) \ce \int_X \norm{F_A}^2 d\operatorname{Vol} \text{ .} \]
The key idea is to apply Morse theory to $\phi$.

\begin{enumerate}
\item The function $\phi$ is a function on the infinite-dimensional
  space $\mathcal{A}$ of all connections. This is an affine-linear
  space, hence contractible.
\item The function $\phi$ is invariant under the infinite-dimensional
  symmetry group of all bundle automorphisms $\mathcal{G} =
  \Map(X,G)$, the so-called group of \emph{gauge transformations}.
\item Inside $\mathcal{G}$ we have the subgroup $\mathcal{G}_0 \subset
  \mathcal{G}$ of \emph{based maps} $X \to G$, which is the kernel of
  $\ev \colon \mathcal{G} \to G$, the evaluation at a base point $x_0
  \in X$ given by $\ev(f)=f(x_0)$. That is, $\mathcal{G}_0$ consists
  of all those gauge transformations which are the identity at $x_0$.

  The restricted group $\mathcal{G}_0$ acts freely on $\mathcal{A}$,
  and so we can reduce to a $G$-action on $\mathcal{A} \bigl/
  \mathcal{G}_0$. Moreover, $\mathcal{G}$-equivariant cohomology on
  $\mathcal{A}$ becomes $G$-equivariant cohomology on
  $\mathcal{A}\bigl/\mathcal{G}_0$.
\item We will apply $\mathcal{G}$-equivariant Morse theory to the
  Yang-Mills functional $\phi$ on the space $\mathcal{A}$.
\end{enumerate}

The critical connections for $\phi$ are the those for which the
curvature $F_A$ is covariantly constant. The absolute minimum appears
when $F_A=0$, i.e.\ when $A$ is \emph{flat} (or more generally
\emph{central harmonic}). For higher critical points, $A$ decomposes.

\begin{example}
Let us consider the simplest case, $n=2$. That is, we consider
rank-$2$ bundles, or $U(2)$-bundles, on a Riemann surface $X$. The
determinant line bundle $\det E$ of a rank-$2$ bundle $E$ has degree
$k = c_1(E) = c_1(\Lambda^2 E)$, and $E$ is topologically non-trivial
whenever $k \neq 0$. Let us assume $k=1$; so we are in a different
component of the moduli space than for $k=0$.

At the absolute minimum, $\mathcal{G}$ acts freely. The moduli space
$M_g(2,1)$ is a manifold and contributes $P_{M_g(2,1)}(t)$. At higher
critical points, the bundle is a direct sum of line bundles, $E \cong
L_1 \oplus L_2$, and $\deg L_1 + \deg L_2 = 1$. Assume without loss of
generality that $\deg L_2 > \deg L_1$. Now $\mathcal{G}$ acts with
isotropy subgroup $U(1)$ and contributes
\[ P^{U(1)}_{J(X) \times J(X)} = \frac{(1+t)^{4g}}{1-t^2} \text{ .} \]
\end{example}

What is the contribution of the total space $\mathcal{A}$? We know
that $H^*_\mathcal{G}\bigl(\mathcal{A}\bigr) =
H^*\bigl(B\mathcal{G}\bigr)$, but how do we calculate this?
Following Atiyah and Bott \cite[\S2]{AtBo2} we have:

\begin{enumerate}
\item $B\mathcal{G} = \Map\bigl(X, B\mathcal{G}\bigr)$.
\item For $G = U(1)$, we have $BG = \mathbb{CP}^\infty$. So
      \[ \Map\bigl(X, \mathbb{CP}^\infty\bigr) = \mathbb{Z} \times
\prod\limits_{2g}S^1 \times \mathbb{CP}^\infty \text{ ,} \]
      and
$$P_{B\mathcal{G}}(t) = (1+t)^{2g}\bigl/(1-t^2)~.$$
\item For $G = U(n)$, we have $U(n) \sim U(1) \times
      S^3 \times \dotsb \times S^{2n-1}$, so
      \begin{equation}\label{eq.PUn}
        P_{B\mathcal{G}}(t) = \frac{\prod_{i=1}^n \bigl(1+t^{2i-1}\bigr)^{2g}}
        {\left(\prod_{i=1}^{n-1} \bigl(1-t^{2i}\bigr)\right)\bigl(1-t^{2n}\bigr)} \text{ .}
      \end{equation}
\end{enumerate}
All of these are ``good cases''. We finish with a computation to prove Theorem \ref{thm.pms}.
\begin{equation}\label{eq.highcrit}
  \frac{t^{2g}(1+t)^{2g}}{(1-t^2)(1-t^4)} = \frac{1}{1-t^2}
  \sum_{i=1}^\infty t^{2g+4i} (1+t)^{2g}
\end{equation}
On the right-hand side we recognise the factors $(1-t^2)^{-1} =
P_{\mathbb{CP}^\infty}(t)$ and $(1+t)^{2g} = P_{J(X)}(t)$. We obtain one big equation
\[ \bigl\{\text{minimum}\bigr\} + \bigl\{\text{higher critical points}
   \bigr\} = \bigl\{\text{total space}\bigr\} \text{ ,} \]
where
\begin{eqnarray*}
  \text{minimum} &=& P_{M^0_g(2,1)} \text{ , the series of the space of interest,}\\
  \text{higher points} &=& \text{the expression \eqref{eq.highcrit}, and}\\
  \text{total space} &=& (1+t^3)^{2g}\bigl/(1-t^2)(1-t^4) \text{ from Equation \eqref{eq.PUn} with $n=2$,}
\end{eqnarray*}
for the Yang-Mills functional $\phi$ on the space of all connections
on $U(2)$-bundles with fixed degree $1$. The contribution from the
higher critical points is given by the $L_1 \oplus L_2$ (with fixed
total degree), which is the origin of the Jacobian factor $P_{J(X)}(t)$.

\begin{remark}
For $n \geq 3$, even if we only want to deal with the co-prime case
$\gcd(n,k)=1$, the inductive step will need a general case (e.g.\
$n=3$, $k=1$ can decompose into $E_2 \oplus E_1$ with $\rk E_i = i$
and $\deg E_2 = 0$, $\deg E_1 = 1$). But Morse theory still works to
give induction if we use equivariant cohomology and equivariant
Poincar\'e series. (The Poincar\'e series $P_M(t)$ will not be a
polynomial).
\end{remark}

\section{Counting rational points}

\subsection{Finite fields}

Fields with finitely many elements are either the integers modulo some
prime $p$, written $\mathbb{F}_p \ce \mathbb{Z}\bigl/p\mathbb{Z}$, or
some algebraic extension thereof, written $\mathbb{F}_q$ with $q=p^n$
for some $n\geq1$. Note that every field is a vector space over its
prime subfield $\mathbb{F}_p$, and the characteristic is in each case
the prime $p$. We can consider an algebraic variety $V$ defined over
any field, in particular over $\mathbb{F}_q$ -- for example by
considering as the defining equations of $V$ polynomials with integer
coefficients and reducing modulo $p$.

\begin{example}[Projective spaces]
Let $V \ce \mathbb{P}\bigl(\mathbb{F}_q^n\bigr) =
\bigl(\mathbb{F}_q^n\setminus\{0\}\bigr) \bigl/
\mathbb{F}_q^\times$. The number of points in $V$ is
\[ N_q(V) = \frac{q^n-1}{q-1} = 1+q+q^2+\dotsb+q^{n-1} \text{ .} \]
Observe:
\begin{enumerate}
\item Over the field $\mathbb{F}_{q^m}$, the number of points is
      \[ N_{q^m}(V) = 1 + q^m + q^{2m} + \dotsb + q^{m(n-1)} \text{ ,} \]
      so varying $m$ determines a polynomial in $q$ via $m \mapsto
      N_{q^m}(V) \in \mathbb{Z}[q]$.
\item Setting $q=t^2$ gives the Poincar\'e polynomial of
      $\mathbb{P}(\mathbb{C}^n) = \mathbb{CP}^{n-1}$. This indicates a
      relation between counting rational points over finite fields and
      Betti numbers of complex varieties.
\item Replacing $q$ by $q^{-1}$ gives
      \[ N_q(V) = \frac{q^n(1-q^{-n})}{q(1-q^{-1})} = q^{n-1}\bigl(1 +
         q^{-1} + \dotsb q^{-(n-1)}\bigr) \text{ ,} \]
      and
      \[ \frac{N_q(V)}{q^{n-1}} = 1 + q^{-1} + q^{-2} + \dotsb +
         q^{-(n-1)} \text{ (Poincar\'e Duality).} \]
\item Let $n\to\infty$. We get $1\bigl/\bigl(1-q^{-1}\bigr)$, and
      putting $q=t^{-2}$ we get $1\bigl/\bigl(1-t^2\bigr) =
      P_{\mathbb{CP}^\infty}(t)$.
\end{enumerate}
\end{example}

\paragraph{Zeta functions.} The $\zeta$-function of an algebraic
variety $V$ over $\mathbb{F}_q$ is
\[ Z_V(t) = \exp\left( \sum_{m=1}^\infty N_{q^m}(V) \frac{t^m}{m}
   \right) \text{ ,} \]
where $N_{q^m}(V)$ is the number of points of $V$ over the finite
field $\mathbb{F}_{q^m}$. We define further
\[ \zeta_V(s) \ce Z_V(q^{-s}) \text{ ,} \]
which is the analogue of the Riemann $\zeta$-function. Note that
$\abs{q^{-s}} = q^{-\Re(s)}$. In the special case where $V$ is a
single point, $Z_V(t) = \frac{1}{1-t}$.

\subsection{The Weil conjectures}

(The Weil conjectures were proved by A.\ Grothendieck and P.\ Deligne.)

\begin{theorem}\label{thm.weil}
Let $V$ be a non-singular projective algebraic variety over a finite
field $\mathbb{F}_q$. Then
\begin{enumerate}
\item $Z_V(t)$ is a rational function of $t$.
\item If $n = \dim V$, then
      \[ Z_V(t) = \frac{p_1(t) \, p_3(t) \, \dotsm \, p_{2n-1}(t)}
         {p_0(t) \, p_2(t) \, \dotsb \, p_{2n}(t)} \text{ ,} \]
      where each root $\omega$ of $p_i$ has $\abs{\omega}=q^{-i/2}$.
\item The roots of $p_i$ are interchanged with the roots of $p_{2n-i}$
      under the substitution $t \to 1\bigl/q^n\,t$.
\item If $V$ is the reduction of an algebraic variety over a subfield
      of $\mathbb{C}$, then the Betti numbers $b_i$ of the variety
      $V(\mathbb{C})$ are $b_i=\deg p_i$.
\end{enumerate}
\end{theorem}

\begin{remark}
Part (2) of Theorem \ref{thm.weil} is the Riemann hypothesis for
function fields. Part (3) is the functional equation for $\zeta(s)$.
\end{remark}

\paragraph{Steps in the proof.}
\begin{enumerate}
\item Define cohomology groups $H^i\bigl(V\bigr)$ which are the
      analogues to $H^i\bigl(V(\mathbb{C})\bigr)$. (Done by Grothendieck.)
\item Use the Frobenius map $\phi \colon V \to V$, $x \mapsto
      x^q$. This maps preserves both multiplication and addition. The
      fixed points of $\phi^m$ are the points of
      $V(\mathbb{F}_{q^m})$, and there are $N_{q^m}(V)$ of them.
\item Apply the Lefschetz fixed point theorem: The number of fixed
      points of a map $f \colon X \to X$ is
      \[ \sum_{i=0}^{\dim X} (-1)^i \tr\bigl(f^* \colon H^i(X;\mathbb{Z})
         \to H^i(X;\mathbb{Z})\bigr) \text{ .} \]
      Take $X=V$, $f=\phi$ and $H^i$ to be Grothendieck cohomology:
      \[ N_{q^m}(V) = \sum_i \tr\bigl((\phi^m)^* \colon H^i(V)
         \to H^i(V) \bigr) = \sum_i (-1)^i \sum_j \omega_{ij}^m \text{ ,} \]
      where the $\omega_{ij}$ are the eigenvalues of $\phi_*$ acting on $H_i(V)$.
\item Now compute:
      \begin{multline*}
        Z_V(t) = \exp\left(\sum_{m=1}^\infty N_{q^m}(V) \frac{t^m}{m}\right) \\
        = \exp\Bigl(\sum_{i} (-1)^i \sum_{j} -\log(1-\omega_{ij}\,t)\Bigr)
        = \prod_{i\text{ odd}} p_i(t) \bigl/ \prod_{i\text{ even}} p_i(t) \text{ ,}
      \end{multline*}
      where $p_i(t) = \prod_j (1-\omega_{ij}\,t)$. This proves the theorem subject to
\item Poincar\'e duality, and
\item the Riemann hypothesis: $\abs{\omega_{ij}}=q^{i/2}$ for all $i,j$ (done by Deligne).
\end{enumerate}

\begin{example}
Let $V=X_g$ be an algebraic curve of genus $g$. Then
\[ Z_V(t) = \frac{\prod_{j=1}^{2g} \bigl(1-\omega_j t\bigr)}
   {\bigl(1-t\bigr) \bigl(1-qt\bigr)} \]
and
\[ \zeta_V(s) = \frac{\prod_{j=1}^{2g} \bigl(1-\omega_j q^{-s}\bigr)}
   {\bigl(1-q^{-s}\bigr)\bigl(1-q^{-s+1}\bigr)} \text{ .} \]
\end{example}

\begin{example}
Let $V = M_g(n,k)$ with $\gcd(n,k)=1$ be the moduli of stable vector
bundles over $X_g$ of rank $n$ and degree $k$. If we can compute
$N_{q^m}(V)$ for all $m$, then Theorem \ref{thm.weil} gives the Betti
numbers of $V(\mathbb{C})$, i.e.\ the Poincar\'e polynomial of
$M_g(n,k)$ over $\mathbb{C}$.
\end{example}

How do we compute the number of points of $M_g(n,k)$ over
$\mathbb{F}_q$? We use two key ideas:
\begin{enumerate}
\item All bundles are trivial if we allow poles (of all orders), i.e.\
      if we work with the field of rational functions on $X_g$.
\item The vector space $A$ of power series over $\mathbb{F}_q$ of the form
      \begin{equation}\label{eq.powerF}
        \sum_{j=0}^\infty a_j t^j \in \mathbb{F}_q \lsem t \rsem
      \end{equation}
      is infinite-dimensional but \emph{compact}, since it is a
      product of finite (hence compact) sets.
\end{enumerate}

The space $A$ has a natural measure $\mu$, which is normalised such
that $\mu(A)=1$. Let $A_r \leq A$ be the linear subspace of power
series of the form \eqref{eq.powerF} which satisfy $a_0 = a_1 = \dotsb
= a_{r-1} = 0$. Then the quotient space $A \bigl/ A_r$ has $q^r$
points, so $\mu(A_r) = q^{-r}$.

We define the infinite projective space over $\mathbb{F}_q$ to be
\[ \mathbb{P}(\mathbb{F}_q^\infty) \equiv \mathbb{F}_q\mathbb{P}^\infty
   \ce (A \setminus\{0\}) \bigl/ \mathbb{F}_q^\times \text{ .} \]
Since $\{0\}$ has measure zero,
\[ \mu\bigl(\mathbb{F}_qP^\infty\bigr) = \mu(A) \bigl/
   \abs{\mathbb{F}_q^\times} = \frac{1}{q-1} \text{ .} \]
(Compare this with the Poincar\'e series $P_{\mathbb{CP}^\infty}(t) =
\frac{1}{1-t}$.)

The way in which we just dealt with infinite dimensions and computed
measures is our inspiration for counting points in moduli spaces over
a finite field $\mathbb{F}_q$: Allowing poles and using measures we
can compute the number of points as ratios of measures.

\begin{example}
The group of isomorphism classes of line bundles over $X_g$ is
isomorphic to the divisor class group $\Cl(X_g)$ of $X_g$, which is
\[ \Cl(X_g) \ce \Div(X_g) \bigl/ \bigl( D \sim D+(f) \bigr) \text{ .} \]
A divisor $D$ is a formal finite sum $D = \sum_j k_j Q_j$, where the
$Q_j \in X_g$ are points and $k_j \in \mathbb{Z}$. Now pick a local
coordinate $u$ near a point $Q$ and let $f$ be a local power series
\[ f(u) = \sum_{k=-N}^\infty a_k \, u^k \text{ , with } a_{-N}\neq 0 \text{ .} \]
Multiplication by elements of a compact group $\mathcal{K}_Q$ reduces
this to $f(u) = u^{-N}$. (The group is the group of holomorphic power
series around $Q$ with non-vanishing constant term, i.e.\ the
invertible elements.) So the group $\Div(X_g)$ of all divisors on
$X_g(\mathbb{F}_q)$ is
\[ \Div(X_g) = \prod_{x \in X_g} \mathcal{K}_x\backslash \mathcal{A}_x
   = \mathcal{K} \backslash \mathcal{A} \text{ ,} \]
and the group of divisor classes of degree $0$, written $\Cl^0(X_g)$, is
\[ \Cl^0(X_g) = \mathcal{K} \backslash \mathcal{A} \slash K^\times \text{ ,} \]
where $K=K(X_g)$ is the function field of $X_g$. The measure
$\mu\bigl(\mathcal{A}\slash K^\times\bigr)$ is finite, and counting
points gives the answer $q^{2g}$.
\end{example}

\paragraph{Bundles of higher rank.} To study the moduli space $M_g(n,k)$ for $n>1$,
i.e.\ the moduli space of bundles of higher rank, we can use the same method, provided
we fix the determinant. We have
\[ \frac{1}{\mu(\mathcal{K})} = (q-1) \sum_E \frac{1}{\abs{\Aut(E)}} \text{ ,} \]
where $\mu$ is the Tamagawa measure (with $c=1$). We have further
\[ \frac{1}{\mu(\mathcal{K})} = q^{(n^2-1)(g-1)} \zeta_{X_g}(2)
   \dotsm \zeta_{X_g}(n) \text{ .} \]
In particular, for $n=2$ and $k=1$ the sum over all bundles $E$ splits
into a sum over stable bundles and a sum over unstable bundles, where
for a stable bundle $E$ we have $\Aut(E) = \{1\}$. Therefore
\[ \sum_E \frac{1}{\abs{\Aut(E)}} = \abs{M^0_g(2,1)} + \sum_{r=1}^\infty \frac{1}{\abs{\Aut(E)}} \text{ ,} \]
where the last sum is a geometric series running over all bundles $E =
L^r \oplus L^{1-r}$ and extensions. This gives an explicit formula for
$\abs{M^0_g(2,1)}$, and hence by the Weil conjectures for $P_{M_g(2,1)}(t)$.

\paragraph{Computing measures.} Let $\alpha$ run over all points of
$M_g^0(n,k)$, i.e.\ orbits of $\mathcal{K}$ acting on
$\mathcal{A}^*\bigl/K^\times$. Then
\[ \sum_\alpha \mu\bigl(\mathcal{K} \bigl/ K_\alpha \bigr)
   = \mu\bigl(\mathcal{A}^*\bigl/K^\times\bigr) = C \text{ ,} \]
or
\[ \sum_{\alpha} \frac{1}{\abs{\mathcal{K}_\alpha}} = \frac{C}{\mu(\mathcal{K})} \text{ .} \]
The only automorphisms of line bundles are scalars, so
$\abs{\mathcal{K}_\alpha} = q-1$. Also,
\[ \sum_{\alpha} \frac{1}{\abs{\mathcal{K}_\alpha}} = \frac{\abs{J(X_g)}}{q-1} \text{ .} \]
We need to know the value of $C$ and $\mu(\mathcal{K})$. Both depend
on the precise normalisation of $\mu$. If we choose $C=1$, then we get
$1\bigl/\mu(\mathcal{K}) = \abs{J(X_g)}$.

\section{Comparison of equivariant Morse theory and counting rational points}

We obtain the same formula for $P_{M}(t)$ and agreement term by term
in the method of the proof. This also works for all $n,k$ and other
groups than $U(n)$. The key points are the following:

\begin{itemize}
\item The total space is ``trivial'': The space of connections is
      affine-linear, hence contractible, and the Tamagawa measure of
      $SL(n;\mathbb{C})$ is $1$.
\item Let $I$ be the isotropy group. We can compute $P_{BI}(t)$ and divide by $\mu(I)$.
\end{itemize}

With $\mathcal{G} = \Map\bigl(X_g, U(n)\bigr) \cong \mathcal{K}$,
\[ P_{B\mathcal{G}}(t) = \prod_{k=1}^n (1+t^{2k-1})^{2g} \Bigl/
   (1-t^{2n}) \prod_{k=1}^{n-1} (1-t^{2k})^2 \text{ ,} \]
and
\[ \frac{1}{\mu(\mathcal{K})} = q^{(n^2-1)(g-1)} \zeta_{X_g}(2) \dotsm \zeta_{X_g}(n) \text{ .} \]
These agree using the formula
\[ \zeta_{X_g}(s) = \prod_{i=1}^{2g}(1 - \omega_i q^{-s})
   \bigl/ (1-q^{-s})(1-q^{1-s}) \text{ .} \]

\paragraph{Questions.}\begin{enumerate}
\item Why do these two formulae agree? (``Quantum analogue of the Weil conjectures'')
\item Is there an extension of the Weil conjectures to infinite
      dimensions?
\item Is computing measures on ad\`elic spaces analogous to Feynman
      integration in gauge theories?
\end{enumerate}

\section{Relation to physics}

Does physics help us understand the questions we raised in the last
section? Is there a relation to the original $\zeta$-function? (This
leads to arithmetic algebraic geometry (Arakelov theory) and further
speculations.)

\paragraph{The Yang-Mills functional} came from physics over
$4$-dimensional space-time. It can be considered formally over a
compact Riemannian manifold $X$ of any dimension $d$. In particular,
\begin{itemize}
\item if $d=2$, $X$ is a Riemann surface and we have many results about moduli spaces;
\item if $d=4$ we have \emph{Donaldson theory}.
\end{itemize}

\paragraph{Quantum field theory.}\begin{enumerate}
\item \emph{Hamiltonian approach:} Consider space and time separately.
      We have a Hilbert space $\mathcal{H}$ of states, and a
      self-adjoint ``Hamiltonian'' operator $H$ acting on
      $\mathcal{H}$. The \emph{evolution} is given by the unitary
      operator $e^{itH}$ on $\mathcal{H}$.
\item \emph{Lagrangian formulation} (relativistically invariant): Let
      $L$ be a functional on some space of functions $f$ on space-time,
      e.g.\ $L(f) = \int \abs{\nabla f}^2$.

\item The \emph{Feynman integral} is $\int \exp \bigl( \frac{i}{\hbar}
      L(f) \bigr)$, integrated over all functions $f$ on $\mathbb{R}^3
      \times [0,\tau]$ with $f(0) = u$ and $f(\tau) = v$, determines
      the value $\left\langle u, e^{i\tau H}v \right\rangle$. This
      relates to the Hamiltonian approach. (Recall that the Lagrangian
      and Hamiltonian are related via the Legendre transform.)
\end{enumerate}

\paragraph{Topological quantum field theories.} For some special
Lagrangians, we get $H=0$, and so time evolution is just the
identity. In this case, the Feynman integrals give topological
information, and we call these cases \emph{topological quantum field
theories}. There are many interesting examples of topological QFTs in
dimensions $2$, $3$ and $4$.

In four dimensions, we get Donaldson theory and Seiberg-Witten theory,
but these have no parameters.

In three dimensions, we get Chern-Simons theory, which does have an
interesting parameter. Let $A$ be a $G$-connection over $X$, where
$G=U(n)$. Let
\[ L = CS(A) = \frac{2\pi}k \int_X \tr\bigl(A \wedge dA +
   \frac23 A \wedge A \wedge A \bigr) \text{ .} \]

The Hilbert space is the space of holomorphic sections of a line
bundle $L^k$ over $M_n(X_g)$, where $X_g$ is a Riemann surface. (This
three-dimensional theory is related to two-dimensional conformal field
theory.) We get topological invariants of $3$-manifolds and knots
inside them (Jones, Witten).

In two dimensions, there is also a Yang-Mills theory with Lagrangian
$L(A) = \norm{F_A}^2 = \int_{X_g} \abs{F_A}^2$. (This is the function
on the space of connections to which we applied equivariant Morse
theory.) This theory is physical and not just topological, but we can
solve it exactly. A coupling constant $\epsilon$ is introduced and the Feynman
integral is formally
\[ Z(\epsilon) = \frac{1}{\operatorname{vol}(\mathcal{G})} \int_{\mathcal{A}}
   \exp\Bigl(-\frac{1}{2\epsilon} \norm{F_A}^2\Bigr) \, dA \text{ .} \]
This has a non-trivial dependence on $\epsilon$ and can be used to
compute the multiplicative structure on the cohomology ring
$H^*\bigl(M(X_g,n)\bigr)$ (Witten).

This quantum field theory looks promising, but does not give the
Poincar\'e series of $M(X_g,n)$. Question: Is there an analogue over a
finite field (where the Frobenius map is related to scaling
$\epsilon$)? Another possibility is to use a (super-symmetric) variant
of Chern-Simons theory for a $3$-manifold $S^1 \times X_g$ (or more
generally a circle bundle or a Seifert fibration). The Hilbert space
is $\Omega^*\bigl(M(X_g,n)\bigr)$, the space of all differential forms
on $M(X_g,n)$, which comes equipped with a differential $d$ and its
adjoint (with respect to the symplectic structure) $d^*$. However,
this seems to involve integration for functions on $S^1 \times X_g$,
while we want just functions on $X_g$ (for the analogy with finite
fields).

A possible idea is contained in Witten-Beasley for another theory of
Chern-Simons type, where integration is reduced to $X_g \subset S^1
\times X_g$ as the fixed-point set of a symmetry.

\section{Finite-dimensional approximations}

We can use approximations to link topology with finite fields and then
pass to a limit. Let us consider approximations to $BG$.

For $G=U(n)$,
\[ BG = \varinjlim_{N\to\infty} \frac{U(N)}{U(n)\times U(N-n)}
   = \varinjlim_{N\to\infty} \Gr_n(\mathbb{C}^N) = \Gr_n(\mathbb{C}^\infty) \text{ .} \]
For maps $f \colon X_g \to BG$, fix a degree $\deg(f) = m$ (and then
let $m \to \infty$). For fixed $N$, $m$, the space of holomorphic maps
$f \colon X_g \to \Gr_n(\mathbb{C}^N)$ of degree $m$ forms a finite-dimensional
algebraic variety $V(N,m)$.

The idea of finite-dimensional approximations is the following:
Holomorphic maps are determined by their behaviour at ``poles'', and
the Gra\ss{}mannians $Gr_n$ can be embedded in projective space. We
can study whether continuous maps can be approximated by holomorphic
maps, apply the Weil conjectures to $V(N,m)$ and take limits.

This is a reasonable programme.

\section{Relation of \texorpdfstring{$\zeta$}{zeta}-functions for
         finite fields and Riemann's \texorpdfstring{$\zeta$}{zeta}-function}

The original Riemann $\zeta$-function is
\[ \zeta(s) = \sum_{n=1}^\infty \frac{1}{n^s} = \sum_{p\text{ prime}}
   \Bigl(1-\frac{1}{p^s}\Bigr)^{-1} \text{ ,} \]
where the last expression is also known as the \emph{Euler product},
whose factors are so-called \emph{local factors}. (They are called
thus with reference to the closed points $(p)$ of the scheme
$\spec(\mathbb{Z})$.) The $\zeta$-function ostensibly contains
information about the set of primes.

Now let $V$ be an algebraic variety over a finite field
$\mathbb{F}_p$.  We want to define a $\zeta$-function for
$V$. If $V = \{*\} = \spec(\mathbb{F}_p)$ is a single point, let
\[ \zeta_V(s) \ce \bigl(1-p^{-s}\bigr)^{-1} \text{ .} \]
In general, if $V$ is any variety defined over $\mathbb{Z}$, we define
\[ \zeta_V(s) \ce \prod_p \zeta_{V_p}(s) \text{ ,} \]
where $V_p$ is the reduction of $V$ modulo $p$. We need to look out
for special ``bad'' primes and add a term for the ``infinite prime''
(arising in valuation theory).

By the Weil conjectures, $\zeta_{V_p}(s)$ is given by a rational
function of $t=p^{-s}$ in terms of the Frobenius action on cohomology.

\begin{example}
Let $V$ be an elliptic curve (i.e.\ of genus $1$) defined over
$\mathbb{Z}$. The Weil formula gives
\[ Z(t) = \frac{(1-\alpha t)(1-\beta t)}{(1-t)(1-pt)} \text{ ,} \]
where $\alpha, \beta$ are eigenvalues of the Frobenius map $\phi$ on
$H^1$, and further we have $\abs{\alpha} = \abs{\beta} = p^{-1/2}$,
$\beta = \alpha^{-1}$ and $\alpha + \beta = a = \tr\bigl(
\phi^*\rvert_{H^1(V)} \bigr)$. Put
\[ L_p(s) = \Bigl(\text{numerator of $Z(t)$ with $t=p^{-s}$}\Bigr)
   = 1-a_p p^{-s} + p^{1-2s} \text{ ,} \]
and
\[ L_V(s) = c.\prod_p L_p(s) \text{ .} \]
\end{example}

\begin{theorem}[Hasse-Weil Conjecture]
With $V$ as above and with suitable choices for the infinite prime and
for bad primes, the function $L_p(s)$ extends holomorphically to all
$s \in \mathbb{C}$, and $L_V(s) = \pm L_V(2s)$.
\end{theorem}

The Hasse-Weil Conjecture has now been proved by Wiles, Taylor and
others. Similar conjectures exist for all $V$ and all $H^i$. (There is
one $L$-function for each $i$.)

\begin{remark}[The ad\`elic picture for $\mathbb{Q}$ or number fields]
This is a comment on the double coset space $\mathcal{K}\backslash G_A
\slash G_K$ used for an algebraic curve over $\mathbb{F}_p$. For
$\mathbb{Q}$ or $\mathbb{Z}$ and for $SL(2)$ we have $SO(2;\mathbb{R})
\backslash SL(2;\mathbb{R})$, which is the upper-half plane (or
hyperbolic plane). The double coset space is
\[ \mathcal{M} \ce = SO(2;\mathbb{R}) \backslash SL(2;\mathbb{R}) \slash SL(2;\mathbb{Z}) \text{ ,} \]
the moduli space of elliptic curves. To compute the area of
$\mathcal{M}$. we start with $SL(2;\mathbb{R}) \slash
SL(2;\mathbb{Z})$, which is a three-dimensional manifold with an
invariant volume form. We decompose it into $SO(2;\mathbb{R})$-orbits
and integrate.
\end{remark}

\section{Arithmetic algebraic geometry (Arakelov theory)}

Suppose we have an algebraic variety $V$ of dimension $d$ defined over
the integers $\mathbb{Z}$. We can either embed $\mathbb{Z}$ into
$\mathbb{C}$ and consider $V(\mathbb{C})$ as a complex variety, or we
can form the residues $\mathbb{Z} \to \mathbb{Z}\bigl/p$ and get a
corresponding variety $V_p$. So in fact we get a family $V_p$ over the
primes in $\mathbb{Z}$, and we include $V_\infty$ sitting over the
infinite prime. This family, a scheme over $\spec \mathbb{Z}$, is an
algebraic variety of dimension $d+1$. If $d=0$ we get a number field,
if $d=1$ we get a so-called \emph{arithmetic surface}.

``In the big picture, physics is at infinity, and number theory at the
finite points.''

We may try to extend theorems from surfaces to their arithmetic
analogues.

\paragraph{Non-Abelian theories.} For $d=0$ and $G=SL(n)$, we
have the Langlands programme, also known as non-Abelian class-field
theory. For $d=1$ we study the local theory at $p$. For $p=\infty$, we
have the \emph{geometric} Langlands programme, which has been related
by Witten to quantum field theories over $V(\mathbb{C})$. What is the
ultimate goal? Perhaps quantum field theories over arithmetic
varieties? One would start with the case $d=1$.

\section{Other questions}

Can we extend our results from curves to varieties of higher
dimensions? Recall that for a curve $X_g$ and gauge group $G=SU(n)$,
we know the Poincar\'e series
\[ P_{\Map(X_g, BG)}(t) = \prod_{k=1}^n (1+t^{2k-1})^{2g} \Bigl/
   (1-t^{2n}) \prod_{k=1}^{n-1} (1-t^{2k})^2 \text{ .} \]
Over $\mathbb{F}_q$,
\[ \operatorname{vol}(\mathcal{K})^{-1} = q^{(n^2-1)(g-1)}
   \zeta_{X_g}(2) \dotsm \zeta_{X_g}(n) \text{ ,} \]
where
\[ \zeta_{X_g}(s) = \prod_{i=1}^{2g}(1 - \omega_i q^{-s})
   \bigl/ (1-q^{-s})(1-q^{1-s}) \]
and $\mathcal{K}$ is the maximal compact subgroup of $G_{A_X}$. Both
formulae extend from curves to varieties $V$ of all dimensions and
still appear to be closely related. We may study, for example, bundles over
\begin{itemize}
\item $\mathbb{P}^2$,
\item $(\mathbb{P}^2, \mathbb{P}^1)$,
\item $\mathbb{P}^1 \times X_g$ (here Morse theory is trickier),
\item $X_g$ with gauge group $\Omega(G)$, this is related to the previous point,
\item and also Weil theory for some infinite-dimensional cases.
\end{itemize}

\end{document}